\documentclass[a4paper,12pt]{article}
\usepackage{amsfonts,amssymb,latexsym,amsmath}
\usepackage{multicol}
\usepackage[cp1251]{inputenc}
\usepackage[russian]{babel}
\begin{document}
\renewcommand{\refname}{References}

\date{}
\title{Invariants of the adjoint action\\ in the nilradical of a
parabolic subalgebra\\ of types $B_n$, $C_n$, $D_n$}
\author{V. Sevostyanova\thanks{This research was partially supported by the RFBR
(project 08-01-00151-а)}}

\maketitle

\begin{center}
\parbox[b]{330pt}{\small\textsc{Abstract.} We study invariants of
the adjoint action of the unipotent group in the nilradical of a
parabolic subalgebra of types $B_n$, $C_n$, $D_n$. We introduce the
notion of expanded base in the set of positive roots and construct
an invariant for every root of the expanded base. We prove that
these invariants are algebraically independent. We also give an
estimate of the transcendence degree of the invariant field.}
\end{center}

\vspace{0.5cm}

Let $K$ be an algebraically closed field of characteristic zero. For
a fixed positive integer~$n$, let $G$ denote one of the following
classical groups defined over~$K$: the symplectic
group~$\mathrm{Sp}_{2n}(K)$, the even orthogonal
group~$\mathrm{O}_{2n}(K)$, and the odd orthogonal
group~$\mathrm{O}_{2n+1}(K)$. Throughout the paper, we set
$$m=\left\{\begin{array}{ll}
2n&\mbox{ if either }G=\mathrm{Sp}_{2n}(K)\mbox{ or
}\mathrm{O}_{2n}(K);\\
2n+1&\mbox{ if }G=\mathrm{O}_{2n+1}(K).\\
\end{array}\right.$$
By $\mathrm{U}_{m}(K)$ denote the upper unitriangular group
consisting of all uppertriangular matrices with unit elements on the
diagonal. Let $\mathrm{B}_m(K)$ be the Borel group consisting of all
upper triangular matrices with nonzero elements on the diagonal.
Denote
$$N=G\cap \mathrm{U}_{m}(K)\mbox{ and }B=G\cap \mathrm{B}_{m}(K).$$

Let $P\supset B$ be a parabolic subgroup of the group~$G$. Denote by
$\mathfrak{p}$, $\mathfrak{b}$, $\mathfrak{n}$ the Lie subalgebras
in the Lie algebra $\mathfrak{g}=\mathrm{Lie}\,G$ that correspond to
the Lie subgroups $P$,~$B$, and~$N$. We represent
\mbox{$\mathfrak{p}=\mathfrak{r}\oplus\mathfrak{m}$} as the direct
sum of the nilradical~$\mathfrak{m}$ and the block diagonal
subalgebra~$\mathfrak{r}$ with size of blocks $(r_1,\ldots,r_s)$.
Consider the adjoint action of the group~$P$ in the
nilradical~$\mathfrak{m}$:
$$\mathrm{Ad}_gx=gxg^{-1},\ x\in\mathfrak{m},\ g\in P.$$
The subalgebra $\mathfrak{m}$ is invariant relative to the adjoint
action of the group~$P$. We extend this action to the regular
representation in the algebra~$K[\mathfrak{m}]$ and in the
field~$K(\mathfrak{m})$:
$$\mathrm{Ad}_g f(x)=f(g^{-1}xg),\ f(x)\in K(\mathfrak{m}),\ g\in P.$$

Since the subalgebra $\mathfrak{m}$ is invariant relative to the
adjoint action of the group~$P$, if follows that $\mathfrak{m}$ is
invariant relative to the action of the Lie subgroup $N$. The
question concerning the structure of the algebra of invariants
$K[\mathfrak{m}]^N$ is open and seems to be a considerable challenge
(see~\cite{S2}). We construct a system of invariants
$$\{M_{\xi},\ \xi\in S,\ L_{\varphi},\ \varphi\in\Phi\}$$
(see Notation~11 and~(\ref{L_q})) in the present paper. We show that
this system of invariants is algebraically independent over~$K$. We
also get an estimate of the transcendence degree of the invariant
field:
$$\mathrm{trdeg}\,K(\mathfrak{m})^N\geqslant|S|+|\Phi|.$$ We show in
paper~\cite{S1} that the estimate for the case of type~$A_n$ is
sharp. This question is open for the other types.

\medskip
Let $T$ be the maximal torus of~$G$ consisting of all diagonal
matrices and $\Delta=\Delta(G,T)$ be the root system defined by~$T$
(see \cite{GG}). By definition, $\Delta$~is a subset of the Abelian
group $X(T)=\mathrm{Hom}(T,K)$ consisting of homomorphisms from~$T$
to~$K$. Let $1\leqslant i\leqslant n$. Denote by $\varepsilon_i$ an
element of the group $X(T)$ such that $\varepsilon_i(t)=t_i$ for all
$t\in T$. Here we denote by $t_i\in K$ the $(i,i)$th entry of the
matrix~$t\in T$. Then $\Delta=\Delta^{\!+}\cup\Delta^{\!-}$ and
$$\Delta^{\!+}=\left\{\begin{array}{ll}
\{\varepsilon_i\pm \varepsilon_j:\ 1\leqslant i<j\leqslant
n\}\cup\{2\varepsilon_i:\ 1\leqslant i\leqslant n\}&\mbox{ if
}G=\mathrm{Sp}_{2n}(K);\\
\{\varepsilon_i\pm \varepsilon_j:\ 1\leqslant i<j\leqslant
n\}&\mbox{ if }G=\mathrm{O}_{2n}(K);\\
\{\varepsilon_i\pm \varepsilon_j:\ 1\leqslant i<j\leqslant
n\}\cup\{\varepsilon_i:\ 1\leqslant i\leqslant n\}&\mbox{ if
}G=\mathrm{O}_{2n+1}(K).\\
\end{array}\right.$$
The roots in $\Delta^{\!+}$ are said to be \emph{positive} (and the
roots in $\Delta^{\!-}$ are said to be \emph{negative}). The system
of positive roots $\Delta^{\!+}_\mathfrak{r}$ of the reductive
subalgebra $\mathfrak{r}$ is a subsystem in~$\Delta^{\!+}$.

We consider the mirror order~$\prec$ on the set $\{0,\pm1,\ldots,\pm
n\}$, which is defined as
$$1\prec2\prec\ldots \prec n\prec0\prec-n\prec\ldots\prec-2\prec-1.$$
We index the rows (from left to right) and columns (from top to
bottom) of any $m\times m$ matrix according to this order.

For any~$\gamma\in\Delta^{\!+}$, we set
$$\mathcal{E}(\gamma)=\left\{\begin{array}{ll}
(-j,-i)&\mbox{if }\gamma=\varepsilon_i-\varepsilon_j,\
1\leqslant i<j\leqslant n;\\
(j,-i)&\mbox{if }\gamma=\varepsilon_i+\varepsilon_j,\
1\leqslant i<j\leqslant n;\\
(i,-i)&\mbox{if }G=\mathrm{Sp}_{2n}(K)\mbox{ and
}\gamma=2\varepsilon_i,\ 1\leqslant i\leqslant n;\\
(0,-i)&\mbox{if }G=\mathrm{O}_{2n+1}(K)\mbox{ and
}\gamma=\varepsilon_i,\
1\leqslant i\leqslant n.\\
\end{array}\right.$$

We define a relation in $\Delta^{\!+}$ for which
$$\gamma'\succ\gamma\mbox{ whenever }\gamma'-\gamma\in\Delta^{\!+}.$$
If $\gamma'\succ\gamma$ or $\gamma'\prec\gamma$, then the
roots~$\gamma$ and $\gamma'$ are said to be \emph{comparable}.

Denote by~$M$ the set of roots $\gamma\in\Delta^{\!+}$ such that the
corresponding subalgebras $\mathfrak{g}_{\gamma}$ are
in~$\mathfrak{m}$. We identify the algebra $K[\mathfrak{m}]$ with
the polynomial algebra in the variables $x_{i,j}$, $i\prec j$,
whenever $(-j,-i)=\mathcal{E}(\gamma)$ and $\gamma\in M$.

\medskip
\textbf{Definition 1.}~A subset $S$ in~$M$ is called a \emph{base}
if the elements in $S$ are not pairwise comparable and for any
$\gamma\in M\setminus S$ there exist $\xi\in S$ such that
$\gamma\succ\xi$.

\medskip
\textbf{Definition 2.}~Let $A$ be a subset in~$S$. We say that
$\gamma$ is a \emph{minimal element} in~$A$ if there is no $\xi\in
A$ such that $\gamma\succ\xi$.

\medskip
Note that $M$ has a unique base~$S$, which can be constructed in the
following way. We form a set~$S_1$ of minimal elements in~$M$. By
definition, $S_1\subset S$. Then we form a set $M_1$, which is
obtained from~$M$ by deleting $S_1$ and all elements
$$\{\gamma\in M:~\exists~\xi\in S_1,~\gamma\succ\xi\}.$$
The subset of minimal elements $S_2$ in~$M_1$ is also contained
in~$S$, and so on. Continuing this process, we get the base~$S$ as a
union of the sets $S_1,S_2,\ldots$

\medskip
\textbf{Lemma 3.} \emph{Let} $\gamma\in M\backslash S$.
\begin{enumerate}
\item \emph{Suppose that} $\gamma=\varepsilon_i+\varepsilon_j$,
$i<j$, \emph{is such that $\mathcal{E}(\xi)$ equals one of the
following values}: $(j,-k)$, \emph{where} $-i\succ -k$, $(k,-i)$,
\emph{where} $j\prec k$, \emph{or} $(k,-j)$ \emph{for some} $k$.
\item \emph{If} $\gamma=\varepsilon_i-\varepsilon_j$, $i<j$,
\emph{then there exists a root $\xi\in S$ such that
$\mathcal{E}(\xi)$ equals} $(-j,k)$, \emph{where} $-i\succ k$,
\emph{or} $(k,-i)$, \emph{where} $j\prec k$.
\item \emph{Let} $\gamma=2\varepsilon_i$, \emph{then there exists
$\xi\in S$ such that} $\mathcal{E}(\xi)=(k,-i)$, \emph{where}
$k\succ i$.
\end{enumerate}

\medskip
\textsc{Proof}. We prove the first item of the lemma. The other
items are proved similarly.

Let $\gamma=\varepsilon_i+\varepsilon_j\in M\backslash S$, $i<j$. By
Definition 1, there exists a root $\xi\in S$ such that
$\gamma-\xi\in\Delta^{\!+}$. We write all suitable roots~$\xi$.
\begin{enumerate}
\item $\xi=\varepsilon_k+\varepsilon_j$, where $i<k$; hence if
$k<j$, then $\mathcal{E}(\xi)=(j,-k)$, if $j<k$, then
$\mathcal{E}(\xi)=(k,-j)$.
\item $\xi=\varepsilon_i+\varepsilon_k$, where $k>j$; then
$\mathcal{E}(\xi)=(k,-i)$.
\item $\xi=\varepsilon_i-\varepsilon_k$, where $i<k$; then
$\mathcal{E}(\xi)=(-k,-i)$.
\item $\xi=\varepsilon_j-\varepsilon_k$, where $j<k$; then
$\mathcal{E}(\xi)=(-k,-j)$.
\item $\xi=\varepsilon_i$; then
$\mathcal{E}(\xi)=(0,-i)$ if $G=\mathrm{O}_{2n+1}(K)$.
\item $\xi=\varepsilon_j$; then
$\mathcal{E}(\xi)=(0,-j)$ if $G=\mathrm{O}_{2n+1}(K)$.
\item $\xi=2\varepsilon_i$; then
$\mathcal{E}(\xi)=(i,-i)$ if $G=\mathrm{Sp}_{2n}(K)$.
\item $\xi=2\varepsilon_j$; then
$\mathcal{E}(\xi)=(j,-j)$ if $G=\mathrm{Sp}_{2n}(K)$.
\end{enumerate}
Thus we have proved the lemma.~$\Box$

\medskip
\textbf{Corollary.} \emph{Let} $G=\mathrm{Sp}_{2n}(K)$, \emph{and
let a number $i>0$ be such that there exists a root} $\gamma\in M$,
\emph{where $\mathcal{E}(\gamma)=(j,-i)$ for some}~$j$. \emph{Then
there exists a root $\xi\in S$ such that $\mathcal{E}(\xi)=(k,-i)$
for some}~$k$.

\medskip
\textsc{Proof}. Suppose that there is no root $\xi\in S$ such that
$\mathcal{E}(\gamma)=(j,-i)$. Consider a root $2\varepsilon_i\in M$.
By Lemma 4, there exists $\xi\in S$ such that
$\mathcal{E}(\xi)=(k,-i)$ for some~$k$, which contradicts the
assumption.~$\Box$

\medskip
\textbf{Definition 4.} An ordered collection of positive roots
$\{\gamma_1,\ldots,\gamma_s\}$ is called a \emph{chain} if
$\mathcal{E}(\gamma_1)=(a_1,a_2)$,
$\mathcal{E}(\gamma_2)=(a_2,a_3)$,
$\mathcal{E}(\gamma_3)=(a_3,a_4)$, and so on.

\medskip
Let $(r_1,\ldots,r_{p-1},r_p,r_{p-1},\ldots,r_1)$ be the sizes of
the blocks in the reductive subalgebra~$\mathfrak{r}$. We denote
$R=\displaystyle\sum_{t=1}^{p-1}r_t$. Let $\gamma\in S$ be a root
such that $\mathcal{E}(\gamma)=(a,-b)$, $b>0$. If $R\prec a\prec-R$,
then $\gamma$ is called a root lying \emph{to the right of the
central block} in the subalgebra~$\mathfrak{r}$.

\medskip
Let $k$ be the greatest number such that the root
$\gamma=\varepsilon_k+\varepsilon_{k+1}$ lies in the root
system~$M$. Denote
$$\Gamma_{\!\mathfrak{r}}=\left\{\begin{array}{ll}
\Delta^{\!+}_{\mathfrak{r}}\cup\{2\varepsilon_i:k<i\leqslant
n\}&\mbox{if }G=\mathrm{O}_{m}(K),\mathrm{O}_{2n}(K);\\
\Delta^{\!+}_{\mathfrak{r}}\backslash\{2\varepsilon_i:k<i\leqslant
n\}&\mbox{if }G=\mathrm{Sp}_{m}(K).\\
\end{array}\right.$$

\medskip
\textbf{Definition 5.}~Assume that one of two roots $\xi,\xi'\in S$
does not lie to the right of the central block in~$\mathfrak{r}$. We
say that two roots $\xi,\xi'$ form an \emph{admissible pair}
$q=(\xi,\xi')$ if there exists a root~$\alpha_q$ in the set
$\Delta_{\mathfrak{r}}^{\!+}$ such that the collection of roots
$\{\xi,\alpha_q,\xi'\}$ is a chain.

Suppose that two roots $\xi,\xi'\in S$ are to the right of the
central block in~$\mathfrak{r}$ and $\mathcal{E}(\xi)=(a,-b)$,
$\mathcal{E}(\xi')=(a',-b')$. Similarly, a pair $q=(\xi,\xi')$ is
called an \emph{admissible pair} if there exists a root
$\alpha_q\in\Gamma_{\!\mathfrak{r}}$ such that
$\mathcal{E}(\alpha_q)=(-a,a')$.

Note that the root $\alpha_q$ can be found from~$q$ uniquely.

\medskip
We form the set $Q=Q(\mathfrak{p})$ that consists of admissible
pairs of roots in~$S$. Let roots~$\xi$ and~$\xi'$ form an admissible
pair. Assume that the roots~$\xi$ и~$\xi'$ are such that
$\mathcal{E}(\xi)=(a,-b)$, $\mathcal{E}(\xi')=(a',-b')$, and
$a\preccurlyeq a'$. For every admissible pair $q=(\xi,\xi')$ we
construct a positive root $\varphi_q=\alpha_q+\xi'$. Consider the
subset $\Phi=\{\varphi_q:~q\in Q\}$.

\medskip
\textbf{Definition 6.} The set $S\cup\Phi$ is called an
\emph{expanded base}.

\medskip
\textbf{Example 7.} In the Lie algebra
$\mathfrak{g}=\mathfrak{o}_{16}(K)$, let the reductive
subalgebra~$\mathfrak{r}$ have the following sizes of diagonal
blocks: $(3,1,2,4,2,1,3)$. Now we write all roots in the expanded
base:
$$S=\big\{\xi_1=\varepsilon_6-\varepsilon_7,\
\xi_2=\varepsilon_5-\varepsilon_8,\
\xi_3=\varepsilon_4-\varepsilon_5,$$
$$\xi_4=\varepsilon_3-\varepsilon_4,\
\xi_5=\varepsilon_2-\varepsilon_6,\
\xi_6=\varepsilon_1+\varepsilon_8\big\}.$$ We write the set of
admissible pairs and the corresponding roots of the system~$\Phi$.
$$Q=\left\{(\xi_1,\xi_1),(\xi_1,\xi_2),(\xi_1,\xi_6),(\xi_2,\xi_2),
(\xi_1,\xi_3)\right\},$$
$$\Phi=\big\{\varphi_{(\xi_1,\xi_1)}=\varepsilon_6+\varepsilon_7,\
\varphi_{(\xi_1,\xi_2)}=\varepsilon_6+\varepsilon_8,\
\varphi_{(\xi_1,\xi_6)}=\varepsilon_6-\varepsilon_8,$$$$
\varphi_{(\xi_21,\xi_2)}=\varepsilon_5+\varepsilon_8,\
\varphi_{(\xi_1,\xi_3)}=\varepsilon_4-\varepsilon_6\big\}.$$

\medskip
From the given parabolic subalgebra, we construct a square diagram
and mark in the diagram each root of the expanded base. Suppose that
the positive root~$\gamma$ corresponds to the pair of integers
$\mathcal{E}(\gamma)=(j,-i)$, $i>0$. Thus we mark the root~$\gamma$
in~$(j,-i)$ entry of the diagram. We label a root of the set~$S$ by
the symbol~$\otimes$ and a root of~$\Phi$ by the symbol~$\times$.
The other entries in the diagram are empty.

\medskip
\textbf{Example 8.} Let $G=\mathrm{O}_{16}(K)$. Let the sizes of
diagonal blocks of the reductive subalgebra be as in Example~7. Then
we have the following diagram:
\medskip
\begin{center}
{\begin{tabular}{|p{0.1cm}|p{0.1cm}|p{0.1cm}|p{0.1cm}|p{0.1cm}|
p{0.1cm}|p{0.1cm}|p{0.1cm}|p{0.1cm}|p{0.1cm}|p{0.1cm}|p{0.1cm}|
p{0.1cm}|p{0.1cm}|p{0.1cm}|p{0.1cm}|c}
\multicolumn{2}{l}{{\footnotesize 1\quad
2}}&\multicolumn{2}{l}{{\footnotesize 3\quad
4}}&\multicolumn{2}{l}{{\footnotesize 5\quad 6}}&
\multicolumn{2}{l}{{\footnotesize 7\quad
8}}&\multicolumn{2}{c}{{\hspace{-10pt}{\footnotesize-8}\
{\footnotesize-7}}}&
\multicolumn{2}{c}{{\hspace{-10pt}{\footnotesize-6}\
{\footnotesize-5}}}&
\multicolumn{2}{c}{{\hspace{-10pt}{\footnotesize-4}\
{\footnotesize-3}}}&
\multicolumn{2}{c}{{\hspace{-10pt}{\footnotesize-2}\ {\footnotesize-1}}}\\
\cline{1-16} \multicolumn{3}{|l|}{1}&\multicolumn{13}{|c|}{}&{\footnotesize 1}\\
\cline{16-16} \multicolumn{3}{|c|}{1}&\multicolumn{12}{|c|}{}&&{\footnotesize 2}\\
\cline{15-16} \multicolumn{3}{|r|}{\qquad\ \ \,1}&\multicolumn{11}{|c|}{}&&&{\footnotesize 3}\\
\cline{1-4}\cline{14-16} \multicolumn{3}{|c|}{}&1&\multicolumn{9}{|c|}{}&&&&{\footnotesize 4}\\
\cline{4-6}\cline{13-16}
\multicolumn{4}{|c|}{}&\multicolumn{2}{|l|}{1}&\multicolumn{6}{|c|}{}&&&&&{\footnotesize 5}\\
\cline{12-16}
\multicolumn{4}{|c|}{}&\multicolumn{2}{|c|}{\quad\,1}&\multicolumn{5}{|c|}{}&&&&&&{\footnotesize 6}\\
\cline{5-16} \multicolumn{6}{|c|}{}&\multicolumn{4}{|l|}{1}&$\times$&&&&&&{\footnotesize 7}\\
\cline{11-16}
\multicolumn{6}{|c|}{}&\multicolumn{4}{|l|}{\quad\ 1}&$\times$&$\times$&&&&$\otimes$&{\footnotesize 8}\\
\cline{11-16} \multicolumn{6}{|c|}{}&\multicolumn{4}{r|}{1\hspace{13pt}\ }&$\times$&$\otimes$&&&&&{\footnotesize -8}\\
\cline{11-16}
\multicolumn{6}{|c|}{}&\multicolumn{4}{|r|}{\qquad\qquad1}&$\otimes$&
&&&&&{\footnotesize -7}\\
\cline{7-16}
\multicolumn{10}{|c|}{}&\multicolumn{2}{l|}{1}&$\times$&&
$\otimes$&&{\footnotesize -6}\\
\cline{13-16} \multicolumn{10}{|c|}{}&\multicolumn{2}{r|}{1}&
$\otimes$&&&&{\footnotesize -5}\\
\cline{11-16} \multicolumn{12}{|c|}{}&1&$\otimes$&&&
{\footnotesize -4}\\
\cline{13-16} \multicolumn{13}{|c|}{}&\multicolumn{3}{|l|}{1}&{\footnotesize -3}\\
\multicolumn{13}{|c|}{}&\multicolumn{3}{|c|}{1}&{\footnotesize -2}\\
\multicolumn{13}{|c|}{}&\multicolumn{3}{|r|}{1}&{\footnotesize -1}\\
\cline{1-16} \multicolumn{16}{c}{Diagram 1}\\
\end{tabular}}
\end{center}

\medskip
\textbf{Example 9.} Let $G=\mathrm{Sp}_{16}(K)$. Let the sizes of
diagonal blocks of the reductive subalgebra be as in Example~8. Then
we get the following diagram.
\medskip
\begin{center}
{\begin{tabular}{|p{0.1cm}|p{0.1cm}|p{0.1cm}|p{0.1cm}|p{0.1cm}|
p{0.1cm}|p{0.1cm}|p{0.1cm}|p{0.1cm}|p{0.1cm}|p{0.1cm}|p{0.1cm}|
p{0.1cm}|p{0.1cm}|p{0.1cm}|p{0.1cm}|c}
\multicolumn{2}{l}{{\footnotesize 1\quad
2}}&\multicolumn{2}{l}{{\footnotesize 3\quad
4}}&\multicolumn{2}{l}{{\footnotesize 5\quad 6}}&
\multicolumn{2}{l}{{\footnotesize 7\quad
8}}&\multicolumn{2}{c}{{\hspace{-10pt}{\footnotesize-8}\
{\footnotesize-7}}}&
\multicolumn{2}{c}{{\hspace{-10pt}{\footnotesize-6}\
{\footnotesize-5}}}&
\multicolumn{2}{c}{{\hspace{-10pt}{\footnotesize-4}\
{\footnotesize-3}}}&
\multicolumn{2}{c}{{\hspace{-10pt}{\footnotesize-2}\ {\footnotesize-1}}}\\
\cline{1-16} \multicolumn{3}{|l|}{1}&\multicolumn{12}{|c|}{}&&{\footnotesize 1}\\
\cline{15-16} \multicolumn{3}{|c|}{1}&\multicolumn{11}{|c|}{}&&&{\footnotesize 2}\\
\cline{14-16} \multicolumn{3}{|r|}{\qquad\ \ \,1}&\multicolumn{10}{|c|}{}&&&&{\footnotesize 3}\\
\cline{1-4}\cline{13-16} \multicolumn{3}{|c|}{}&1&\multicolumn{8}{|c|}{}&&&&&{\footnotesize 4}\\
\cline{4-6}\cline{12-16}
\multicolumn{4}{|c|}{}&\multicolumn{2}{|l|}{1}&\multicolumn{5}{|c|}{}&&&&&&{\footnotesize 5}\\
\cline{11-16}
\multicolumn{4}{|c|}{}&\multicolumn{2}{|c|}{\quad\,1}&\multicolumn{4}{|c|}{}&&&&&&&{\footnotesize 6}\\
\cline{5-16} \multicolumn{6}{|c|}{}&\multicolumn{4}{|l|}{1}&&&&&&&{\footnotesize 7}\\
\cline{11-16}
\multicolumn{6}{|c|}{}&\multicolumn{4}{|l|}{\quad\ 1}&$\times$&&&&&$\otimes$&{\footnotesize 8}\\
\cline{11-16} \multicolumn{6}{|c|}{}&\multicolumn{4}{r|}{1\hspace{13pt}\ }&$\times$&$\otimes$&&&&&{\footnotesize -8}\\
\cline{11-16}
\multicolumn{6}{|c|}{}&\multicolumn{4}{|r|}{\qquad\qquad1}&$\otimes$&
&&&&&{\footnotesize -7}\\
\cline{7-16}
\multicolumn{10}{|c|}{}&\multicolumn{2}{l|}{1}&$\times$&&
$\otimes$&&{\footnotesize -6}\\
\cline{13-16} \multicolumn{10}{|c|}{}&\multicolumn{2}{r|}{1}&
$\otimes$&&&&{\footnotesize -5}\\
\cline{11-16} \multicolumn{12}{|c|}{}&1&$\otimes$&&&
{\footnotesize -4}\\
\cline{13-16} \multicolumn{13}{|c|}{}&\multicolumn{3}{|l|}{1}&{\footnotesize -3}\\
\multicolumn{13}{|c|}{}&\multicolumn{3}{|c|}{1}&{\footnotesize -2}\\
\multicolumn{13}{|c|}{}&\multicolumn{3}{|r|}{1}&{\footnotesize -1}\\
\cline{1-16} \multicolumn{16}{c}{Diagram 2}\\
\end{tabular}}
\end{center}

\medskip
\textbf{Пример 10.} Let $G=\mathrm{O}_{17}(K)$. Let
$(1,2,5,1,5,2,1)$ be the sizes of diagonal blocks in the reductive
subalgebra. Then we have the following diagram.
\medskip
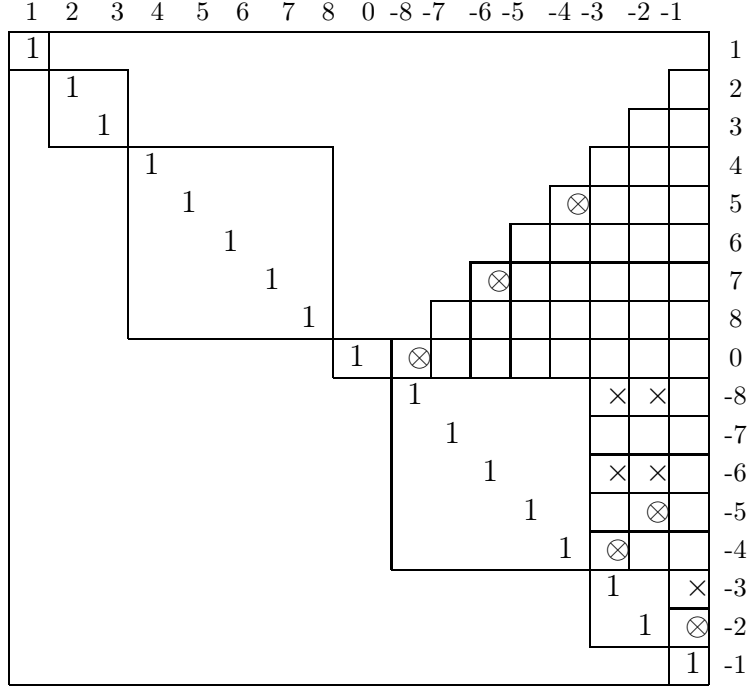
\begin{center}
{\begin{tabular}{|p{0.1cm}|p{0.1cm}|p{0.1cm}|p{0.1cm}|p{0.1cm}|
p{0.1cm}|p{0.1cm}|p{0.1cm}|p{0.1cm}|p{0.1cm}|p{0.1cm}|p{0.1cm}|p{0.1cm}|
p{0.1cm}|p{0.1cm}|p{0.1cm}|p{0.1cm}|c}
\multicolumn{2}{l}{{\footnotesize 1\quad
2}}&\multicolumn{2}{l}{{\footnotesize 3\quad
4}}&\multicolumn{2}{l}{{\footnotesize 5\quad 6}}&
\multicolumn{3}{l}{{\footnotesize 7\quad 8\quad
0}}&\multicolumn{2}{c}{{\hspace{-10pt}{\footnotesize-8}\
{\footnotesize-7}}}&
\multicolumn{2}{c}{{\hspace{-10pt}{\footnotesize-6}\
{\footnotesize-5}}}&
\multicolumn{2}{c}{{\hspace{-10pt}{\footnotesize-4}\
{\footnotesize-3}}}&
\multicolumn{2}{c}{{\hspace{-10pt}{\footnotesize-2}\ {\footnotesize-1}}}\\
\cline{1-17} 1&\multicolumn{16}{|c|}{}&{\footnotesize 1}\\
\cline{1-3}\cline{17-17} &\multicolumn{2}{|l|}{1}&\multicolumn{13}{|c|}{}&&{\footnotesize 2}\\
\cline{16-17} &\multicolumn{2}{|r|}{\quad 1}&\multicolumn{12}{|c|}{}&&&{\footnotesize 3}\\
\cline{2-8} \cline{15-17} \multicolumn{3}{|c|}{}&\multicolumn{5}{|l|}{1}&\multicolumn{6}{|c|}{}&&&&{\footnotesize 4}\\
\cline{14-17} \multicolumn{3}{|c|}{}&\multicolumn{5}{|l|}{\hspace{10pt} 1}&\multicolumn{5}{|c|}{}&$\otimes$&&&&{\footnotesize 5}\\
\cline{13-17} \multicolumn{3}{|c|}{}&\multicolumn{5}{|c|}{1}&\multicolumn{4}{|c|}{}&&&&&&{\footnotesize 6}\\
\cline{12-17} \multicolumn{3}{|c|}{}&\multicolumn{5}{|r|}{1 \hspace{10pt}\ }&\multicolumn{3}{|c|}{}&$\otimes$&&&&&&{\footnotesize 7}\\
\cline{11-17} \multicolumn{3}{|c|}{}&\multicolumn{5}{|r|}{\hspace{55pt} 1}&\multicolumn{2}{|c|}{}&&&&&&&&{\footnotesize 8}\\
\cline{4-17} \multicolumn{8}{|c|}{}&1&$\otimes$&&&&&&&&{\footnotesize 0}\\
\cline{9-17} \multicolumn{9}{|c|}{}&\multicolumn{5}{|l|}{1}&$\times$&$\times$&&{\footnotesize -8}\\
\cline{15-17} \multicolumn{9}{|c|}{}&\multicolumn{5}{|l|}{\hspace{10pt} 1}&&&&{\footnotesize -7}\\
\cline{15-17} \multicolumn{9}{|c|}{}&\multicolumn{5}{|c|}{1}&$\times$&$\times$&&{\footnotesize -6}\\
\cline{15-17} \multicolumn{9}{|c|}{}&\multicolumn{5}{|r|}{1\hspace{13pt}\ }&&$\otimes$&&{\footnotesize -5}\\
\cline{15-17} \multicolumn{9}{|c|}{}&\multicolumn{5}{|r|}{1}&$\otimes$&&&{\footnotesize -4}\\
\cline{10-17} \multicolumn{14}{|c|}{}&\multicolumn{2}{|l|}{1}&$\times$&{\footnotesize -3}\\
\cline{17-17} \multicolumn{14}{|c|}{}&\multicolumn{2}{|r|}{1}&$\otimes$&{\footnotesize -2}\\
\cline{15-17} \multicolumn{16}{|c|}{}&1&{\footnotesize -1}\\
\cline{1-17} \multicolumn{16}{c}{Diagram 3}\\
\end{tabular}}
\end{center}

\medskip
We construct the formal matrix $\mathbb{X}$ as follows. Let
$\gamma\in M$. Let $\mathcal{E}(\gamma)=(-j,-i)$, where $i>0$. Then
the variable $x_{i,j}$ occupies the position $(i,j)$. The position
$(-j,-i)$ is occupied by the variable $x_{i,j}$, where
$G=\mathrm{Sp}_{2n}(K)$ and $j<0$, or by the variable $-x_{i,j}$ in
the other cases. If $\mathcal{E}(\gamma)=(i,-i)$, $i>0$, then the
variable $x_{i,-i}$ stands in the position $(i,-i)$.

\medskip
Assume that $\gamma\in M$ and $\mathcal{E}(\gamma)=(a,-b)$,
$a\in\mathbb{Z}$, $b>0$. Denote by $S_{\gamma}$ the set of $\xi$
in~$S$ such that $\mathcal{E}(\xi)=(i,j)$, $i\succ a$, and $j\prec
-b$. Let $S_{\gamma}=\{\xi_1,\ldots,\xi_k\}$,
$\mathcal{E}(\xi_i)=(a_i,-b_i)$, where $a_i\in\mathbb{Z}$, $b_i>0$
for any $1\leqslant i\leqslant k$.

Let $b_{i_1},b_{i_2},\ldots,b_{i_s}$ be numbers in the set
$\{b_1,b_2,\ldots,b_k\}$ that are greater than the number~$a$.

\medskip
\textbf{Notation 11}. Denote by $M_\gamma$ the minor~$M_I^J$ of the
matrix~$\mathbb{X}$ with ordered systems of rows~$I$ and
columns~$J$, where
$$I=\mathrm{ord}\{b_1,\ldots,b_k,b,a_{i_1},\ldots,a_{i_s}\},\quad
J=\mathrm{ord}\{-a,-a_1,\ldots,-a_k,-b_{i_1},\ldots,-b_{i_s}\}.$$
Denote by $\overline{M}_\gamma$ the minor~$M_{I'}^{J'}$ of the
matrix~$\mathbb{X}$, where
$$I'=\mathrm{ord}\{a,a_1,\ldots,a_k,b_{i_1},\ldots,b_{i_s}\},\quad
J'=\mathrm{ord}\{-b_1,\ldots,-b_k,-b,-a_{i_1},\ldots,-a_{i_s}\}.$$
Obviously, $M_{\gamma}=\pm\overline{M}_{\gamma}$.

\medskip
\textbf{Example 12.} Let the group~$G$ and the parabolic subalgebra
be as in Example~10. Let the root $\gamma$ be the root
$\varepsilon_6+\varepsilon_7$; then $S_{\gamma}=\{\varepsilon_8\}$.
We have $\mathcal{E}(\varepsilon_6+\varepsilon_7)=(7,-6)$,
$\mathcal{E}(\varepsilon_8)=(0,-8)$. Since $8\succ7$, we have
$I=\{7,8,0\}$, $J=\{0,-8,-6\}$, whence we get
$$M_{\gamma}=\left|\begin{array}{ccc}
x_{7,0}&x_{7,-8}&x_{7,-6}\\
x_{8,0}&0&-x_{6,-8}\\
0&-x_{8,0}&-x_{6,0}\\
\end{array}\right|.$$

\medskip
Note that the statement $a_i\neq b_j$ is valid for any numbers
$i\neq j$ such that $1\leqslant i,j\leqslant k$. Indeed, we have
$a_i>0$, consequently, if $a_i=b_j$, then
$\xi_i=\varepsilon_{b_i}+\varepsilon_{a_i}$ and
$\xi_j=\varepsilon_{b_j}+\varepsilon_{a_j}$ if $a_j>0$, and
$\xi_j=\varepsilon_{b_j}-\varepsilon_{-a_j}$ if $a_j<0$. Hence
$\xi_i$ and~$\xi_j$ of the base~$S$ are the comparable roots:
$$\xi_j-\xi_i=\varepsilon_{b_i}+\varepsilon_{a_i}-(\varepsilon_{a_i}
\pm\varepsilon_{\pm a_j})=\varepsilon_{b_i}\mp\varepsilon_{\pm
a_j}.$$ We have a contradiction with the definition of a base. Thus
the sets of rows~$I$ and columns~$J$ of the minor $M_{\gamma}$ do
not contain equal elements. Therefore $M_{\gamma}$ and
$\overline{M}_{\gamma}$ are the square minors.

\medskip
Suppose that a set $A$ contains some numbers of the set
$$\{1,2,\ldots,n,0,-n,\ldots,-2,-1\}.$$ Denote
$$\min A=\{k\in A:k\preccurlyeq a\mbox{ for any }a\in A\},$$
$$\max A=\{k\in A:k\succcurlyeq a\mbox{ for any }a\in A\}.$$

The following lemma is needed in the sequel.

\medskip
\textbf{Lemma 13.} \emph{Let} $\xi\in S$, \emph{and rows~$I$ and
columns~$J$ form the minor} $M_{\xi}$. \emph{Let}
$\mathcal{E}(\xi)=(a,-b)$, $b>0$.
\begin{enumerate}
\item \emph{Assume that the number $i\not\in I$ such that} $\min
I\prec i\prec\max I$; \emph{then} $i=a$.
\item \emph{For any number~$j$ such that} $\min J\prec
j\prec\max J$, \emph{we have} $j\in J$.
\end{enumerate}

\textsc{Proof}.
\begin{enumerate}
\item We prove the first statement of the lemma. Let $i\not\in
I$, then, by the definition of the minor~$M_{\xi}$, we conclude that
there is no root $\gamma\in S$ such that
$\mathcal{E}(\gamma)=(j,-i)$ and $i\preccurlyeq a$. Assume that
$i\prec a$.

Consider roots $\gamma_j\in\Delta^{\!+}$ such that
$\mathcal{E}(\gamma_j)=(j,-i)$. Denote by~$A$ the set of the numbers
$j$ such that the roots $\gamma_j$ lie in $M$. As was said above,
for any number $j\in A$ we have $\gamma_j\not\in S$. By the
definition of a base, we deduce that there exists a root $\xi_j\in
S$ such that $\gamma_j-\xi_j\in\Delta^{\!+}$ for the given root
$\gamma_j$. Thus any root $\alpha\in M$ such that
$\mathcal{E}(\alpha)=(j,-k)$, where $j\in A$ and $-k\succ-i$, is
comparable with the root $\xi_j$ of the base~$S$, i.e., any such
root $\alpha$ does not lie in~$S$. Further, since $\min I=b\prec i$,
we have $-i\prec-b$. Since for any number $j\in J$ we have
$j\prec\max J=-a$, it follows that $a\prec\min A$. By the assumption
$i\prec a$, we have $i\prec\min A$. For any root $\alpha\in M$ such
that $\mathcal{E}(\alpha)=(\min A,-i)$ we have $\min A\prec i$. We
obtain a contradiction, and thus $i=a$.

\item Now we prove the second statement. Let~$j$ be a number such
that $\min J\prec j\prec\max J$. Obviously, $\max J=-a$. Consider a
root $\gamma\in M$ such that $\mathcal{E}(\gamma)=(-j,-b)$. Then
$$\gamma=\left\{\begin{array}{ll}
\varepsilon_b-\varepsilon_j&\mbox{if }j>0;\\
\varepsilon_b+\varepsilon_{-j}&\mbox{if }j<0.\\
\end{array}\right.$$
Clearly, $\gamma\not\in S$. Otherwise the root $\gamma$ is
comparable with the root $\xi\in S$. We obtain a contradiction with
the definition of a base.

Let $\gamma=\varepsilon_b-\varepsilon_j$. From Lemma~3 it follows
that there exists a root $\xi'\in S$ such that
$\mathcal{E}(\xi')=(-j,-i)$ for some $-i\prec-b$. Since $j\prec\max
J=-a$, we have $-j\succ a$ and $-i\prec-b$. From the last
inequalities it follows that $\xi'\in S_{\xi}$. Hence, $j\in J$.

Similarly, if $\gamma=\varepsilon_b+\varepsilon_{-j}$, then from
Lemma~3 it follows that there exists a root $\xi'\in S$ such that
either $\mathcal{E}(\xi')=(-j,-i)$ for some $-i\prec-b$ or
$\mathcal{E}(\xi')=(i,j)$ for some $i\succ -j$. In the first case,
we obtain $\xi'\in S_{\xi}$ and $j\in J$. In the second case, we
have $i\succ-j\succ a$ and $j\prec-a\prec-b$. Consequently, $\xi'\in
S_{\xi}$. Since $-j\succ a$, by the definition of the minor
$M_{\xi}$, we have $j\in J$.$~\Box$
\end{enumerate}

\medskip
Suppose a root $\varphi\in\Phi$ corresponds to an admissible pair
$q=(\xi,\xi')\in Q$. We construct a polynomial
\begin{equation}
L_{\varphi}=
\sum_{\scriptstyle\alpha_1,\alpha_2\in\Gamma_{\!\mathfrak{r}}\cup\{0\}
\atop\scriptstyle\alpha_1+\alpha_2=\alpha_q}M_{\xi+\alpha_1}
\overline{M}_{\alpha_2+\xi'}.\label{L_q}
\end{equation}
Theorem 15 shows that the polynomials $L_{\varphi}$ are
$N$-invariants.

\medskip
\textbf{Example 14.} We continue the calculations of Example~7-8. We
construct some polynomials, using the roots of the expanded base.
$$M_{\varepsilon_6-\varepsilon_7}=x_{6,7},\
M_{\varepsilon_5-\varepsilon_8}=\left|\begin{array}{cc}
x_{5,7}&x_{5,8}\\
x_{6,7}&x_{6,8}\\
\end{array}\right|,\
M_{\varepsilon_4-\varepsilon_5}=x_{4,5},\
M_{\varepsilon_3-\varepsilon_4}=x_{3,4},$$
$$M_{\varepsilon_2-\varepsilon_6}=\left|\begin{array}{ccc}
x_{2,4}&x_{2,5}&x_{2,6}\\
x_{3,4}&x_{3,5}&x_{3,6}\\
0&x_{4,5}&x_{4,6}\\
\end{array}\right|,$$
$$L_{\varepsilon_6+\varepsilon_7}=-x_{6,7}x_{6,-7}-x_{6,8}x_{6,-8},$$
$$L_{\varepsilon_6+\varepsilon_8}=-x_{6,-8}\left|\begin{array}{cc}
x_{5,7}&x_{5,8}\\x_{6,7}&x_{6,8}\\
\end{array}\right|-x_{6,8}\left|\begin{array}{cc}
x_{5,7}&x_{5,-8}\\x_{6,7}&x_{6,-8}\\
\end{array}\right|-x_{6,7}\left|\begin{array}{cc}
x_{5,7}&x_{5,-7}\\x_{6,7}&x_{6,-7}\\
\end{array}\right|.$$
$$L_{\varepsilon_5+\varepsilon_8}=2\left|\begin{array}{cc}
x_{5,7}&x_{5,-8}\\x_{6,7}&x_{6,-8}\\
\end{array}\right|\left|\begin{array}{cc}
x_{5,7}&x_{5,8}\\x_{6,7}&x_{6,8}\\
\end{array}\right|,$$
$$L_{\varepsilon_4-\varepsilon_6}=-x_{4,6}x_{6,7}-x_{4,5}x_{5,7}.$$

\medskip
By $E_{i,j}$ denote the standard elementary square matrix having
unit in the $(i,j)$th entry and zeros in all other positions. To
every root $\alpha\in\Delta^+$ corresponds a one-parameter subgroup
$g_{\alpha}(t)$ of square $m\times m$ matrices, where $t\in K$:
\begin{equation}
g_{\alpha}(t)=\left\{\begin{array}{ll}
E+t(E_{i,j}-E_{-j,-i})&\mbox{if }\alpha=\varepsilon_i-\varepsilon_j;\\
E+t(E_{i,-j}+E_{j,-i})&\mbox{if }\alpha=\varepsilon_i+\varepsilon_j\mbox{ and }G=\mathrm{Sp}_{2n}(K);\\
E+t(E_{i,-j}-E_{j,-i})&\mbox{if }\alpha=\varepsilon_i+\varepsilon_j\mbox{ and }G\neq\mathrm{Sp}_{2n}(K);\\
E+tE_{i,-i}&\mbox{if }\alpha=2\varepsilon_i\mbox{ and }G=\mathrm{Sp}_{2n}(K);\\
E+t(E_{i,0}-E_{0,-i})-\frac{t^2}{2}E_{i,-i}&\mbox{if }\alpha=\varepsilon_i\mbox{ and }G=\mathrm{O}_{2n+1}(K).\\
\end{array}\right.\label{g_alpha}
\end{equation}

Now we prove the main statement of the paper.

\medskip
\textbf{Theorem 15.} \emph{For any parabolic subalgebra}, \emph{the
system of polynomials}
$$\{M_\xi,~\xi\in S;~L_{\varphi},~\varphi\in\Phi\}$$
\emph{is contained in $K[\mathfrak{m}]^N$ and is algebraically
independent over}~$K$.

\medskip
\textsc{Proof}. It is well known that for any fixed ordering of
positive roots, any element $g\in N$ can be written in the form
$$g=\prod_{\alpha\in\Delta^{\!+}}g_{\alpha}(t_{\alpha}),$$
where $t_{\alpha}\in K$ such that $\alpha\in\Delta^{\!+}$ are
uniquely determined. Therefore it is sufficient to prove that for
any $\xi\in S$ and $\varphi\in\Phi$ the polynomials $M_\xi$ and
$L_{\varphi}$ are invariants of the adjoint action of the
one-parameter subgroups $g_{\alpha}(t)$ for any
$\alpha\in\Delta^{\!+}$, $t\in K$. Note that the adjoint action by
the element $E+tE_{i,j}\in \mathrm{U}_m(K)$, $i<j$, reduces to the
composition of two transformations: the row~$j$ multiplied by~$-t$
is added to the row~$i$; the column~$i$ multiplied by~$t$ is added
to the column~$j$.

\medskip
Let us show that the minor $M_{\xi}$, $\xi\in S$, is an invariant of
the adjoint action of $g_{\alpha}(t)$. Let
$\mathcal{E}(\xi)=(a,-b)$, $b>0$. Assume that the rows~$I$ and the
columns~$J$ form the minor $M_{\xi}$. We shell give two comments,
the second one is a result of Lemma~13:
\begin{enumerate}
\item The elements $(i,j)$, where $\min I\preccurlyeq
i\preccurlyeq\max I$ and $1\preccurlyeq j\prec\min J$, and $(l,k)$,
where $\max I\prec l\preccurlyeq n$ и $\min J\preccurlyeq
k\preccurlyeq\max J$, of the matrix $\mathbb{X}$ equal zero.
\item Suppose $\min I\preccurlyeq
i\preccurlyeq\max I$; hence if $i\neq a$, then $i\in I$. If $\min
J\preccurlyeq j\preccurlyeq\max J$, then $j\in J$.
\end{enumerate}
If for any number~$i$ such that $\min I\prec i\prec\max I$ we have
$i\in I$, then the above remarks imply that the minor $M_{\xi}$ is
an invariant.

Suppose $i\not\in I$ and $\min I\prec i\prec\max I$; then, by
Lemma~13, we have $i=a$. By the corollary of Lemma~3, we have
$G=\mathrm{O}_m(K)$. Let $I=\{a_1,\ldots,a_k\}$. We prove that the
minor $M_{\xi}$ is an invariant. It is sufficient to show that
$M_{\xi}$ is an invariant of the adjoint action of $g_{\alpha}(t)$,
where the root $\alpha\in\Delta^{\!+}$ is
$\mathcal{E}(\alpha)=(-a,-a_t)$ for some~$t$. Then
$$\mathrm{Ad}_{g_{\alpha}}M_{\xi}=M_{\xi}\pm tM^J_{I'},$$
where $I'=\mathrm{ord}\{a_1,\ldots,a_{t-1},a,a_{t+1},\ldots,a_k\}$.
Let us show that the minor $M^J_{I'}$ equals zero. Consider the sets
$$\widetilde{I}=\mathrm{ord}\{l\in I':~l\succ a\}\subset I'\mbox{
and }\widetilde{J}=\mathrm{ord}\{l:~-l\in\widetilde{I}\}.$$ Since
$a=-\max J$, we get $\widetilde{J}\subset J$. Hence $M^J_{I'}$ is
the minor
\begin{center}
{$M^J_{I'}=$ \begin{tabular}{|c|c|}
\raisebox{3pt}[17pt]{$M_{I'\backslash\widetilde{I}}^{J\backslash\widetilde{J}}$}&
\raisebox{3pt}[17pt]{$M_{I'\backslash\widetilde{I}}^{\widetilde{J}}$}\\
\hline \raisebox{3pt}[17pt]{0}&\raisebox{3pt}[17pt]{$M_{\widetilde{I}}^{\widetilde{J}}$}\\
\end{tabular} .}
\end{center}
Thus,
$M^J_{I'}=M_{I'\backslash\widetilde{I}}^{J\backslash\widetilde{J}}\cdot
M_{\widetilde{I}}^{\widetilde{J}}$. We show that the order of the
minor $M_{\widetilde{I}}^{\widetilde{J}}$ is an odd number. If this
is the case, then $M_{\widetilde{I}}^{\widetilde{J}}$ is an
antisymmetric minor and the order of
$M_{\widetilde{I}}^{\widetilde{J}}$ is an odd number, whence
$M_{\widetilde{I}}^{\widetilde{J}}$ equals zero. Now let $\gamma\in
S$, $\mathcal{E}(\gamma)=(c,-d)$, be a root such that
$d\in\widetilde{I}$. Then $d\succ a$; therefore by the definition
of~$M_{\xi}$ implies that $c\in I$. Since $a\prec d\prec c$, we have
$c\in\widetilde{I}$. Thus, any root $\gamma\in S$, whenever
$\mathcal{E}(\gamma)=(c,-d)$ and $d\in\widetilde{I}$, adds two rows
and two columns to the sets~$I$ и~$J$, respectively. But the
number~$a$ is in the set $\widetilde{I}$; therefor
$M_{\widetilde{I}}^{\widetilde{J}}$ is an antisymmetric minor of odd
order. Hence $M_{\widetilde{I}}^{\widetilde{J}}$ equals zero.
Consequently, $M_{\xi}$ is an invariant for any $\xi\in S$.

\medskip
Let us show that $L_{\varphi}$ is an invariant for any
$\varphi\in\Phi$. Let $q=(\xi,\xi')$ be an admissible pair and
$\{\xi,\alpha_q,\xi'\}$ be a chain. Suppose that
$\mathcal{E}(\xi)=(a,-b)$ and $\mathcal{E}(\xi')=(a',-b')$, where
$b,b'>0$, are such that $a\preccurlyeq a'$. Then
$\varphi=\alpha_q+\xi'$.

\begin{enumerate}
\item
Consider a case where $a'\succcurlyeq-R$. We have $a\neq a'$;
otherwise the roots $\xi$ and $\xi'$ in~$S$ are comparable. Since
$a'<0$, we have $\xi'=\varepsilon_{b'}-\varepsilon_{-a'}$. The roots
$\xi,\alpha_q,\xi'$ form a chain; therefor
$\mathcal{E}(\alpha_q)=(-b,a')$, i.e.,
$\alpha_q=\varepsilon_{-a'}-\varepsilon_b$. Hence
$\varphi=\varepsilon_{b'}-\varepsilon_b$ and
$\mathcal{E}(\varphi)=(-b,-b')$.

Assume that we have the adjoint action of the subgroup
$g_{\alpha}(t)$ on the polynomial $L_{\varphi}$, where
$\mathcal{E}(\alpha)=(-j,-i)$, $i>0$. If $-j\prec-b$ or
$a'\preccurlyeq-i$, then the action of $g_{\alpha}(t)$ does not
change $L_{\varphi}$. Therefore let
$-b\preccurlyeq-j\prec-i\preccurlyeq a'$; then $j>0$ and
$\alpha=\varepsilon_i-\varepsilon_j$. Let the root $\alpha_q$ be the
sum of two roots
$$\alpha_q=\gamma_1+\alpha+\gamma_2,\mbox{ where
}\gamma_1=\varepsilon_j-\varepsilon_b,\
\gamma_2=\varepsilon_{-a'}-\varepsilon_i.$$ Since
$\alpha_q\in\Delta^{\!+}_{\mathfrak{r}}$, we have
$\gamma_1,\gamma_2\in\Delta^{\!+}_{\mathfrak{r}}$.

The following forms can be obtained by direct calculation.
$$\begin{array}{rl}
\mathcal{E}(\xi+\gamma_1+\alpha)=(a,-i),&\mathcal{E}(\gamma_2+\xi')=(-i,-b'),\\
\mathcal{E}(\xi+\gamma_1)=(a,-j),&\mathcal{E}(\alpha+\gamma_2+\xi')=(-j,-b').\\
\end{array}$$
Using (\ref{g_alpha}), we have
\begin{equation}
\begin{array}{rcl}
\mathrm{Ad}_{g_{\alpha}(t)}M_{\xi+\gamma_1+\alpha}&=&M_{\xi+\gamma_1+\alpha}-
tM_{\xi+\gamma_1},\\
\mathrm{Ad}_{g_{\alpha}(t)}\overline{M}_{\alpha+\gamma_2+\xi'}&=&\overline{M}_{\alpha+\gamma_2+\xi'}+
t\overline{M}_{\gamma_2+\xi'}.
\end{array}\label{Ad_g_alpha}
\end{equation}
Applying (\ref{Ad_g_alpha}) to (\ref{L_q}), we obtain
$$\mathrm{Ad}_{g_{\alpha}(t)}L_{\varphi}-L_{\varphi}=
\left(M_{\xi+\gamma_1+\alpha}-
tM_{\xi+\gamma_1}\right)\overline{M}_{\gamma_2+\xi'}+$$$$+
M_{\xi+\gamma_1}\left(\overline{M}_{\alpha+\gamma_2+\xi'}+
t\overline{M}_{\gamma_2+\xi'}\right)-
M_{\xi+\gamma_1+\alpha}\overline{M}_{\gamma_2+\xi'}-
M_{\xi+\gamma_1}\overline{M}_{\alpha+\gamma_2+\xi'}=0,$$ i.e.,
$L_{\varphi}$ is an invariant.

\item Now suppose $a'\prec-R$. Obviously, since
$\mathcal{E}(\alpha_q)=(-a,a')$, we have $a'<0$. Therefore
$\xi'=\varepsilon_{b'}-\varepsilon_{-a'}$. Now we write all
variations of the roots~$\xi$ and~$\alpha_q$:
$$\xi=\left\{\begin{array}{ll}
\varepsilon_b-\varepsilon_{-a}&\mbox{if }a<0;\\
\varepsilon_b+\varepsilon_{a}&\mbox{if }a>0;\\
\varepsilon_b,&\mbox{if }a=0;\\
\end{array}\right.\quad
\alpha_q=\left\{\begin{array}{ll}
\varepsilon_{-a'}+\varepsilon_{-a}&\mbox{if }a<0;\\
\varepsilon_{-a'}-\varepsilon_{a}&\mbox{if }a>0;\\
\varepsilon_{-a'}&\mbox{if }a=0.\\
\end{array}\right.$$

Consider the adjoint action of $g_{\alpha}(t)$. We prove that the
polynomial~$L_{\varphi}$ is an invariant if
$\alpha=\varepsilon_i-\varepsilon_j$ and $\alpha=\varepsilon_i$. The
other cases $\alpha=\varepsilon_i+\varepsilon_j$ and
$\alpha=2\varepsilon_i$ are proved similarly.

\begin{enumerate}
\item Suppose that $\alpha=\varepsilon_i-\varepsilon_j$, $i<j$.
Then $\mathcal{E}(\alpha)=(-j,-i)$. Obviously, $-a\prec a'$. The
action of $g_{\alpha}(t)$ does not change the polynomial
$L_{\varphi}$ if $-i\succcurlyeq a'$ and $-j\prec-a$. Assume that
$-a\preccurlyeq -j\prec-i\prec a'$.

First suppose $a>0$. We write the root $\alpha_q$ as the sum of
three roots, where $\alpha$ is one of these roots and the other
roots are in~$\Gamma_{\!\mathfrak{r}}$.
$$\alpha_q=\varepsilon_{-a'}-\varepsilon_{a}=\gamma_1+\alpha+\gamma_2,$$
where $\gamma_1=\varepsilon_j-\varepsilon_a$ and
$\gamma_2=\varepsilon_{-a'}-\varepsilon_i$. We have
\begin{equation}
\begin{array}{rl}
\xi+\gamma_1+\alpha=\varepsilon_b+\varepsilon_i,&
\gamma_2+\xi'=\varepsilon_{b'}-\varepsilon_i,\\
\xi+\gamma_1=\varepsilon_b+\varepsilon_j,&
\alpha+\gamma_2+\xi'=\varepsilon_{b'}-\varepsilon_j.\\
\end{array}\label{alpha_q=}
\end{equation}
Let $a<0$. As before, the root $\alpha_q$ is the sum of three roots:
$$\alpha_q=\varepsilon_{-a'}+\varepsilon_{-a}=\gamma_1+\alpha+\gamma_2,$$
where $\gamma_1=\varepsilon_{-a}+\varepsilon_j$ and
$\gamma_2=\varepsilon_{-a'}-\varepsilon_i$. Similarly, the
relations~(\ref{alpha_q=}) are verified by direct calculations. If
$a=0$, then $\alpha_q=\varepsilon_{-a'}$ is the sum of roots
$\gamma_1=\varepsilon_j$, $\alpha$ and
$\gamma_2=\varepsilon_{-a'}-\varepsilon_i$. As above, we
get~(\ref{alpha_q=}). Every case yields
relations~(\ref{Ad_g_alpha}). Then
$$\mathrm{Ad}_{g_{\alpha}(t)}L_{\varphi}-L_{\varphi}=0.$$

\item Now suppose that $\alpha=\varepsilon_i$; then,
by~(\ref{g_alpha}), we have
$$g_{\alpha}(t)=E+t\Big(E_{i,0}-E_{0,-i}\Big)-\frac{t^2}{2}E_{i,-i}.$$
The action of $g_{\alpha}(t)$ does not change the polynomial
$L_{\varphi}$ if $a'\preccurlyeq-i$ or $a\geqslant0$. Hence let
$a<0$ and $a'\succ-i$. Then
$\alpha_q=\varepsilon_{-a'}+\varepsilon_{-a}$. We write the
root~$\alpha_q$ as the sum of two roots in $\Gamma_{\!\mathfrak{r}}$
such that a certain root of the sum equals $\alpha+\gamma$,
$\gamma\in\Gamma_{\!\mathfrak{r}}$:
$$\alpha_q=\gamma_1+\gamma_2=\gamma_3+\gamma_4=\gamma_5+\gamma_6,$$
where
$$\begin{array}{lll}
\gamma_1=\varepsilon_{-a}-\varepsilon_{i},&
\gamma_3=\varepsilon_{-a}+\varepsilon_{i},&\gamma_5=\varepsilon_{-a},\\
\gamma_2=\varepsilon_{-a'}+\varepsilon_{i},&
\gamma_4=\varepsilon_{-a'}-\varepsilon_{i},&\gamma_6=\varepsilon_{-a'}.\\
\end{array}$$
Then
$$\begin{array}{lll}\xi+\gamma_1=\varepsilon_b-\varepsilon_i,&
\xi+\gamma_3=\varepsilon_b+\varepsilon_i,&
\xi+\gamma_5=\varepsilon_b,\\
\gamma_2+\xi'=\varepsilon_{b'}+\varepsilon_i,&
\gamma_4+\xi'=\varepsilon_{b'}-\varepsilon_i,&
\gamma_6+\xi'=\varepsilon_{b'}.\\
\end{array}$$
The adjoint action of the subgroup $g_{\alpha}(t)$ on the
corresponding minors gives the following relations:
$$\begin{array}{l}\mathrm{Ad}_{g_{\alpha}(t)}\overline{M}_{\gamma_2+\xi'}=
\mathrm{Ad}_{g_{\alpha}(t)}\overline{M}_{\varepsilon_{b'}+\varepsilon_i}=
\overline{M}_{\varepsilon_{b'}+\varepsilon_i}-t\overline{M}_{\varepsilon_{b'}}-
\frac{t^2}{2}\overline{M}_{\varepsilon_{b'}-\varepsilon_i},\\
\mathrm{Ad}_{g_{\alpha}(t)}M_{\xi+\gamma_3}=
\mathrm{Ad}_{g_{\alpha}(t)}M_{\varepsilon_{b}+\varepsilon_i}=
M_{\varepsilon_{b}+\varepsilon_i}-tM_{\varepsilon_{b}}-
\frac{t^2}{2}M_{\varepsilon_{b}-\varepsilon_i},\\
\mathrm{Ad}_{g_{\alpha}(t)}M_{\xi+\gamma_5}=
\mathrm{Ad}_{g_{\alpha}(t)}M_{\varepsilon_{b}}=
M_{\varepsilon_{b}}+tM_{\varepsilon_{b}-\varepsilon_i},\\
\mathrm{Ad}_{g_{\alpha}(t)}\overline{M}_{\gamma_6+\xi'}=
\mathrm{Ad}_{g_{\alpha}(t)}\overline{M}_{\varepsilon_{b'}}=
\overline{M}_{\varepsilon_{b'}}+t\overline{M}_{\varepsilon_{b'}-\varepsilon_i}.\\
\end{array}$$
Applying these expressions to~(\ref{L_q}), we get
$$\mathrm{Ad}_{g_{\alpha}(t)}L_{\varphi}-L_{\varphi}=
M_{\varepsilon_b-\varepsilon_i}\Big(\overline{M}_{\varepsilon_{b'}+\varepsilon_i}-
t\overline{M}_{\varepsilon_{b'}}-
\frac{t^2}{2}\overline{M}_{\varepsilon_{b'}-\varepsilon_i}\Big)+
$$$$+\Big(M_{\varepsilon_{b}+\varepsilon_i}-tM_{\varepsilon_{b}}-
\frac{t^2}{2}M_{\varepsilon_{b}-\varepsilon_i}\Big)\overline{M}_{\varepsilon_{b'}-\varepsilon_i}+
$$$$+\Big(M_{\varepsilon_{b}}+tM_{\varepsilon_{b}-\varepsilon_i}\Big)
\Big(\overline{M}_{\varepsilon_{b'}}+t\overline{M}_{\varepsilon_{b'}-\varepsilon_i}\Big)-$$$$-
M_{\varepsilon_b-\varepsilon_i}\overline{M}_{\varepsilon_{b'}+\varepsilon_i}-
M_{\varepsilon_{b}+\varepsilon_i}\overline{M}_{\varepsilon_{b'}-\varepsilon_i}-
M_{\varepsilon_{b}}\overline{M}_{\varepsilon_{b'}}=0.$$
\end{enumerate}
\end{enumerate}

To complete the proof, it remains to show that for any $\xi\in S$
and $\varphi\in\Phi$, the polynomials $M_{\xi}$ and $L_{\varphi}$
are algebraically independent. Let
$$E_{\gamma}=\left\{\begin{array}{ll}
E_{i,j}-E_{-j,-i}&\mbox{if }\gamma=\varepsilon_i-\varepsilon_j;\\
E_{i,-j}-E_{j,-i}&\mbox{if }\gamma=\varepsilon_i+\varepsilon_j\mbox{ and }G\neq\mathrm{Sp}_{2n}(K);\\
E_{i,-j}+E_{j,-i}&\mbox{if }\gamma=\varepsilon_i+\varepsilon_j\mbox{ and }G=\mathrm{Sp}_{2n}(K);\\
E_{i,-i}&\mbox{if }\gamma=2\varepsilon_i\mbox{ and }G=\mathrm{Sp}_{2n}(K);\\
E_{i,0}-E_{0,-i}&\mbox{if }\gamma=\varepsilon_i\mbox{ and }G=\mathrm{O}_{2n+1}(K).\\
\end{array}\right.$$
Denote by $\mathcal{Y}$ the subset in~$\mathfrak{m}$ of matrices of
the form
$$\sum_{\xi\in S}c_{\xi}E_{\xi}+\sum_{\varphi\in\Phi}
c'_{\varphi}E_{\varphi}.$$ Consider the restriction homomorphism
$\pi(f)=f|_{\mathcal{Y}}$ of the algebra $K[\mathfrak{m}]$
to~$K[\mathcal{Y}]$. The image $K[\mathfrak{m}]$ of the homomorphism
is the polynomial algebra of $x_{i,j}$, where
$(-j,-i)=\mathcal{E}(\gamma)$ for some root $\gamma\in S\cup\Phi$.
Denote $x_{i,j}=x_{\xi}$ if $\mathcal{E}(\xi)=(-j,-i)$. We have the
following images of the polynomials $M_{\xi}$ and $L_{\varphi}$ for
any $\xi\in S$ and $\varphi\in\Phi$:
$$\pi(M_{\xi})=\pm x_{\xi}\prod_{\gamma\in S_{\xi}}x_{\gamma}^{\delta_{\gamma}},$$
and if a root $\varphi\in\Phi$ corresponds to the admissible pair
$(\xi,\xi')$, then
$$\pi(L_{\varphi})=\pm x_{\varphi}\prod_{\gamma\in\{\xi\}\cup
S_{\xi}\cup S_{\xi'}}x_{\gamma}^{\delta_{\gamma}},$$ where
$\delta_{\gamma}$ equals~1 or~2. Since the system of $\pi(M_{\xi})$
and $\pi(L_{\varphi})$ is algebraically independent, where $\xi\in
S$ and $\varphi\in\Phi$, then the system of $M_{\xi}$ and
$L_{\varphi}$ is algebraically independent.~$\Box$

\medskip
From Theorem 15 we obtain the following consequence.

\medskip
\textbf{Theorem 16.} The dimension of $N$-orbit in~$\mathfrak{m}$ is
no greater than the number
$$\mathrm{dim}\,\mathfrak{m}-|S|-|\Phi|.$$

\textsc{Department of Mechanics and Mathematics, Samara State
University, Russia}\\ \emph{E-mail address}:
\verb"victoria.sevostyanova@gmail.com"

\end{document}